\documentclass[10pt,twoside]{article}
\usepackage{amsmath}%
\usepackage{amsfonts}%
\usepackage{amssymb}%
\usepackage{graphicx}
\usepackage{Latex-document}
\newtheorem{theorem}{Theorem}

\newtheorem{proposition}[theorem]{Proposition}

\title{\bf Localization-Delocalization Phenomena \vskip -2mm  for Random Interfaces\vskip 6mm}
\author{Erwin Bolthausen\vspace*{-0.5cm}\thanks{Institut f\"ur Mathematik, University
of Z\"{u}rich, Winterthurerstrasse 190, 8057 Z\"{U}RICH,
Switzerland. E-mail: eb@amath.unizh.ch}}
\date{\vspace{-8mm}}

\begin{document}

\maketitle

\thispagestyle{first} \setcounter{page}{25}

\begin{abstract}\vskip 3mm

We consider $d$-dimensional random surface models which for $d=1$ are the standard (tied-down) random walks
(considered as a random ``string''). In higher dimensions, the one-dimensional (discrete) time parameter of the
random walk is replaced by the $d$-dimensional lattice $\mathbb{Z}^{d}$, or a finite subset of it. The random
surface is represented by real-valued random variables $\phi_{i},$ where $i\in\mathbb{Z}^{d}.$ A class of natural
generalizations of the standard random walk are gradient models whose laws are
(formally) expressed as%
\[
P\left(  d\phi\right)  =\frac{1}{Z}\exp\left[ -\sum\nolimits_{\left| i-j\right|  =1}V\left(
\phi_{i}-\phi_{j}\right)  \right]  \prod_{i}d\phi_{i},
\]
$V:\mathbb{R\rightarrow R}^{+},$ convex, and with some growth conditions.

Such surfaces have been introduced in theoretical physics as (simplified) models for random interfaces separating
different phases. Of particular interest are localization-delocalization phenomena, for instance for a surface
interacting with a wall by attracting or repulsive interactions, or both together. Another example are so-called
heteropolymers which have a noise-induced interaction.

Recently, there had been developments of new probabilistic tools for such problems. Among them are:

\begin{itemize}
\item Random walk representations of Helffer-Sj\"{o}strand type,
\item Multiscale analysis,
\item Connections with random trapping problems and large deviations.
\end{itemize}

We give a survey of some of these developments.

\vskip 4.5mm


\end{abstract}

\vskip 12mm

\section{Introduction\label{Chap_Intro}}

\vskip-5mm \hspace{5mm}

Gradient models are an important class of random interfaces and random
surfaces. In the mathematical physics literature they are often called
``effective interface models''. The (discrete) random surface is described by
random variables $\left(  \phi_{x}\right)  _{x\in V},$ where $V$ is
$\mathbb{Z}^{d}$ or a subset of it. The $\phi_{x}$ itself are either
$\mathbb{Z}$-valued or $\mathbb{R}$-valued. We will mainly concentrate on the
latter situation which is easier in some respects. If $V$ is a finite subset
of $\mathbb{Z}^{d},$ the law $P_{V}$ of $\phi=\left(  \phi_{x}\right)  _{x\in
V}$ is described via a Hamiltonian%
\begin{equation}
H_{V}\left(  \phi\right)  \overset{\mathrm{def}}{=}\frac{1}{2}\sum_{x,y\in
V}p\left(  y-x\right)  U\left(  \phi_{x}-\phi_{y}\right)  +\sum_{x\in
V,\,y\notin V}p\left(  y-x\right)  U\left(  \phi_{x}\right)
,\label{Hamiltonian}%
\end{equation}
where $U:\mathbb{R\rightarrow R}^{+}$ is symmetric and convex, and
$p:\mathbb{Z}^{d}\rightarrow\left[  0,1\right]  $ is a symmetric probability
distribution on $\mathbb{Z}^{d}.$ The above choice of the Hamiltonian
corresponds to $0$ boundary conditions. Of course, one can consider more
general ones, where the second summand is replaced by $\sum_{x\in V,\,y\notin
V}p\left(  y-x\right)  U\left(  \phi_{x}-\psi_{y}\right)  ,$ $\psi$ being a
configuration outside $V.$ We will be mainly interested in the nearest
neighbor case $p\left(  x\right)  =1/2d,$ for $\left|  x\right|  =1,$ and
$p\left(  x\right)  =0$ otherwise, but more general conditions can also be
considered. We always assume that the matrix $\mathcal{Q}=\left(
q_{ij}\right)  $ given by%
\begin{equation}
q_{ij}\overset{\mathrm{def}}{=}\sum_{x}x_{i}x_{j}p\left(  x\right)
\label{Def_Cov-von-p}%
\end{equation}
is positive definite, and that $p$ has exponentially decaying tails.
Furthermore, the random walk $\left(  \eta_{t}\right)  _{t\in\mathbb{N}}$ with
transition probabilities $p$ is assumed to be irreducible. The Hamiltonian
defines a probability distribution on $\mathbb{R}^{V}$ by%
\begin{equation}
P_{V}\left(  d\phi\right)  \overset{\mathrm{def}}{=}\frac{1}{Z_{V}}\exp\left[
-H_{V}\left(  \phi\right)  \right]  \prod_{x\in V}d\phi_{x},\label{Def_PV}%
\end{equation}
where $d\phi_{x}$ denotes the Lebesgue measure. $Z_{V}$ is the norming
constant%
\begin{equation}
Z_{V}\overset{\mathrm{def}}{=}\int_{\mathbb{R}^{V}}\exp\left[  -H_{V}\left(
\phi\right)  \right]  \prod_{x\in V}d\phi_{x}.\label{Def_Zustandssumme}%
\end{equation}

In the one-dimensional case $d=1,$ $P_{V}$ is the law of a tied down random
walk: Let $\xi_{i},\;i\geq1,$ be i.i.d. random variables with the density
$\operatorname*{const}\times\mathrm{e}^{-U\left(  x\right)  }.$ If $V=\left\{
1,\ldots,n\right\}  ,$ then $P_{V}$ is the law of the sequence $\left(
\sum_{j=1}^{i}\xi_{j}\right)  _{1\leq i\leq n},$ conditioned on $\sum
_{j=1}^{n+1}\xi_{j}=0.$

A special case is the \textit{harmonic }one with $U\left(  x\right)
=x^{2}/2.$ Then $P_{V}$ is a Gaussian measure on $\mathbb{R}^{V}$ which is
centered for $0$-boundary conditions. We usually write $P_{V}^{\mathrm{harm}}$
in this case. The law is therefore given by its covariances%
\[
\gamma_{V}\left(  x,y\right)  \overset{\mathrm{def}}{=}\int\phi_{x}\phi
_{y}\,dP_{V}^{\mathrm{harm}}.
\]
These covariances have a random walk representation: If $V$ is a finite set
then%
\begin{equation}
\gamma_{V}\left(  x,y\right)  =\mathbb{E}_{x}\left(  \sum\nolimits_{s=0}%
^{\tau_{V}-1}1_{y}\left(  \eta_{s}\right)  \right)  ,\label{RWRepresentation}%
\end{equation}
where $\left(  \eta_{s}\right)  _{s\geq0}$ under $\mathbb{P}_{x}$ is a
discrete time random walk on $\mathbb{Z}^{d}$ starting at $x$ and with
transition probabilities $\mathbb{P}_{x}\left(  \eta_{1}=y\right)  =p\left(
y-x\right)  .$ $\tau_{V}$ is the first exit time from $V.$ As a consequence of
this representation one sees that the thermodynamic limit%
\begin{equation}
P_{\infty}^{\mathrm{harm}}\overset{\mathrm{def}}{=}\lim_{n\rightarrow\infty
}P_{V_{n}}^{\mathrm{harm}},\;V_{n}\overset{\mathrm{def}}{=}\left\{
-n,-n+1,\ldots,n\right\}  ^{d}\label{Def_P-infinity}%
\end{equation}
exists for $d\geq3.$ $P_{\infty}^{\mathrm{harm}}$ is the centered Gaussian
measure on $\mathbb{R}^{\mathbb{Z}^{d}}$ whose covariances are given by the
Green's function of the random walk. It is important to notice that this
random field has slowly decaying correlations:%
\[
\gamma_{\infty}\left(  x,y\right)  \approx\frac{\operatorname*{const}}{\left|
x-y\right|  ^{d-2}},\;\left|  x-y\right|  \rightarrow\infty.
\]

For $d=2,$ the thermodynamic limit does not exist, and in fact%
\[
E_{V_{n}}^{\mathrm{harm}}\left(  \phi_{0}^{2}\right)  \approx
\operatorname*{const}\times\log n,\;n\rightarrow\infty.
\]
For $d=1,$ the variance grows of course like $n$ in the bulk. The harmonic
surface is therefore localized for $d\geq3,$ but not for $d=1,2.$

Many of these properties carry over to non-harmonic cases with a convex and
symmetric interaction function $U$ in (\ref{Hamiltonian}). Of particular
importance is that there is a generalization of the representation
(\ref{RWRepresentation})$,$ the Helffer-Sj\"{o}strand representation, see
\cite{HeSj}. The random walk $\left(  \eta_{s}\right)  $ has to be replaced by
a random walk in a dynamically changing random environment. Using this
representation, many of the results for the harmonic case can be generalized
to the case of a convex $U,$ although often not in a quantitatively as precise
form as in the harmonic case. For a probabilistic description of the
Helffer-Sj\"{o}strand representation, see \cite{DeGiIo}.

The main topic of this paper are effects arising from interactions of the
random surface $\left(  \phi_{x}\right)  $ with a ``wall''. The simplest case
of such a wall is the configuration $\phi\equiv0.$ There are many type of
interactions which had been considered in the literature, both in physics and
in mathematics. The simplest one is a local attraction of the surface to this
wall. It turns out that an arbitrary weak attraction localizes the random
field in a strong sense, and in all dimensions. This will be discussed in a
precise way in Section \ref{Chap_Pinning}. Interesting
localization-delocalization phenomena may occur when mixed attractive and
repulsive interactions are present, with phase transitions depending on the
parameters regulating the strength of the interactions. Naturally, these
phenomena are best understood for the one-dimensional case. A simple example
is the following one, which is discussed in details in \cite{Fisher}: Let
$\phi_{0}=0,\phi_{1},\ldots,\phi_{2n-1},\phi_{2n}=0$ be a discrete time
$\mathbb{Z}$-valued, and tied-down, simple random walk, i.e. $P_{n}$ is simply
the uniform distribution on all such paths which satisfy $\left|  \phi
_{x}-\phi_{x-1}\right|  =1.$ Introducing an arbitrary pinning to the wall in
the form%
\[
\hat{P}_{n,\beta}\left(  \phi\right)  =\frac{1}{\hat{Z}_{n,\beta}}\exp\left[
\beta\sum\nolimits_{x=1}^{2n-1}1_{0}\left(  \phi_{x}\right)  \right]
P_{n}\left(  \phi\right)  ,\;\beta>0
\]
strongly localizes the ``random string'', i.e. $\sup_{n,x}\hat{E}_{n,\beta
}\left(  \phi_{x}^{2}\right)  <\infty$ holds for all $\beta>0.$ Furthermore,
the correlations $\hat{E}_{n,\beta}\left(  \phi_{x}\phi_{y}\right)  $ are
exponentially decaying in $\left|  x-y\right|  ,$ uniformly in $n.$ These
facts are easily checked.

On the other hand, if the string is confined to be on one side of the wall,
the situation is completely different. Let $\Omega_{2n}^{+}\overset
{\mathrm{def}}{=}\left\{  \phi:\phi_{x}\geq0,\;1\leq x\leq2n-1\right\}  ,$ and
$\hat{P}_{n,\beta}^{+}\left(  \cdot\right)  \overset{\mathrm{def}}{=}\hat
{P}_{n,\beta}\left(  \cdot\mid\Omega_{2n}^{+}\right)  .$ Then there is a
critical $\beta_{c}>0$ such that the above localization property holds for
$\beta>\beta_{c},$ but not for $\beta<\beta_{c},$ where the path measure
converges, after Brownian rescaling, to the Brownian excursion. For a proof of
this so called ``wetting transition'', see \cite{Fisher}. More precise
information has been obtained recently in this one-dimensional situation in
\cite{Jap}.

There are similar phase transitions for more complicated models. Some of them
will be discussed in Section \ref{Chap_EntropicRep} and Section
\ref{Chap_Copolymer}. We begin in the next section by discussing the pinning
effect alone mainly in the difficult two-dimensional case.

\section{Pinning of two-dimensional gradient fields\label{Chap_Pinning}}

\vskip-5mm \hspace{5mm}

We consider now a gradient field (\ref{Def_PV}), but we modify it by
introducing an attractive local pinning to the wall $\left\{  \phi
\equiv0\right\}  .$ This is often done by modifying the Hamiltonian in the
following way: Let $\psi:\mathbb{R\rightarrow R}^{-}$ be symmetric and with
compact support. Then we put%
\begin{equation}
H_{V,\psi}\left(  \phi\right)  =H_{V}\left(  \phi\right)  +\sum_{x\in V}%
\psi\left(  \phi_{x}\right)  .\label{bump-pinning}%
\end{equation}
Evidently, the corresponding finite volume Gibbs measure favours surfaces
which have the tendency to stick close to the wall. It should be emphasized
that this is a much weaker attraction than in a so-called massive field, where
one takes $\psi$ to be convex, for instance $\psi\left(  x\right)  =x^{2}.$ A
formally slightly easier model can be obtained by not changing the
Hamiltonian, but replacing the Lebesgue measure as the reference measure by a
mixture of the Lebesgue measure and a Dirac measure at $0.$ This corresponds
to the following probability measure on $\mathbb{R}^{V}:$%
\begin{equation}
\hat{P}_{V,\varepsilon}\left(  d\phi\right)  \overset{\mathrm{def}}{=}\frac
{1}{\hat{Z}_{V,\varepsilon}}\exp\left[  -H_{V}\left(  \phi\right)  \right]
\prod_{x\in V}\left(  d\phi_{x}+\varepsilon\delta_{0}\left(  d\phi_{x}\right)
\right)  ,\;\varepsilon>0.\label{Def_PV-pinned}%
\end{equation}
This measure can be obtained from measures defined by the Hamiltonian
(\ref{bump-pinning}) via an appropriate limiting procedure. The nice feature
of (\ref{Def_PV-pinned}) is that $\hat{P}_{V,\varepsilon}$ can trivially be
expanded into a mixture of ``free'' measures: We just have to expand out the
product:%
\begin{align}
\hat{P}_{V,\varepsilon}\left(  d\phi\right)    & =\sum_{A\subset V}%
\varepsilon^{\left|  V\backslash A\right|  }\frac{Z_{A}}{\hat{Z}%
_{V,\varepsilon}}\frac{1}{Z_{A}}\exp\left[  -H_{V}\left(  \phi\right)
\right]  \prod_{x\in A}d\phi_{x}\prod_{x\in V\backslash A}\delta_{0}\left(
d\phi_{x}\right)  \label{Expand_Pinning}\\
& =\sum_{A\subset V}\varepsilon^{\left|  V\backslash A\right|  }\frac{Z_{A}%
}{\hat{Z}_{V,\varepsilon}}P_{A}\left(  d\phi\right)  ,\nonumber
\end{align}
where $P_{A}$ is the measure defined by (\ref{Def_PV}), extended by $0$
outside $A.$ Remark that%
\[
\nu_{V,\varepsilon}\left(  A\right)  \overset{\mathrm{def}}{=}\varepsilon
^{\left|  V\backslash A\right|  }\frac{Z_{A}}{\hat{Z}_{V,\varepsilon}}%
\]
defines a probability distribution on the set of subsets of $V.$ Therefore, we
have represented $\hat{P}_{V,\varepsilon}$ as a mixture of free measures
$P_{A}.$ It should be remarked that similar but technically more involved
expansions are possible also in the case of the Hamiltonian
(\ref{bump-pinning}). The case of $\psi\left(  x\right)  =-a1_{\left[
-b,b\right]  }\left(  x\right)  $ is discussed in \cite{BoVe}. Probably, more
general cases could be handled with the help of the
Brydges-Fr\"{o}hlich-Spencer random walk representation (see \cite{BrFrSp}),
but the results presented here have not been derived in this more general
case. For the sake of simplicity, we stick here to the $\delta$-pinning case
(\ref{Def_PV-pinned}).

The above representation easily leads to a representation of the covariances
of the pinned field. This is particularly simple in the harmonic case
$U\left(  x\right)  =x^{2}/2,$ where one gets%
\[
\int\phi_{x}\phi_{y}\hat{P}_{V,\varepsilon}^{\mathrm{harm}}\left(
d\phi\right)  =\sum_{A\subset V}\nu_{V,\varepsilon}\left(  A\right)
\mathbb{E}_{x}\left(  \sum\nolimits_{s=0}^{\tau_{A}-1}1_{y}\left(  \eta
_{s}\right)  \right)  .
\]
The problem is therefore reduced to a problem of a random walk among random
traps: The distribution $\nu_{V,\varepsilon}$ defines a random trapping
configuration, let's denote it by $\mathcal{A},$ i.e. $P_{\mathrm{trap}%
}\left(  \mathcal{A}=V^{c}\cup\left(  V\backslash A\right)  \right)
\overset{\mathrm{def}}{=}\nu_{V,\varepsilon}\left(  A\right)  ,$ and the
covariances of our pinned measure are given in terms of the discrete Green's
function among these random traps which are killing the random walk when it
enters one of these traps. A difficult point is a precise analysis of the
distribution of $\mathcal{A},$ and a crucial step is a comparison with
Bernoulli measures. The two-dimensional case is the most difficult one. In
three and more dimensions, a comparison of the distribution of $\mathcal{A}$
with a Bernoulli measure is quite easy.

It turns out that the pinning localizes the field in a strong sense. First of
all, the variance of the variables stay bounded as $V\uparrow\mathbb{Z}^{d}.$
Secondly, there is exponential decay of the covariances, uniformly in $V.$
Results of this type have a long history. For $d\geq3,$ and for the harmonic
case with pinning of the type (\ref{bump-pinning}), the localization has been
obtained in \cite{BrFrSp}. In \cite{DuMaRiRo}, boundedness of the absolute
first moment has been proved for $d=2$, but no exponential decay of the
correlations. The first proof of exponential decay of correlations in the
two-dimensional case has been obtained in \cite{BoBr} for the harmonic case.
One drawback of the method used there was that it uses reflection positivity,
which holds only under restrictive assumptions on $p.$ Also, periodic boundary
conditions are required, and so the results are not directly valid for the
$0$-boundary case. A satisfactory approach had then been obtained in
\cite{DeVe} and \cite{IoVe}. The quantitatively precise results presented here
are from \cite{BoVe}, where the critical exponents for the depinning
transition $\varepsilon\rightarrow0$ have been derived, including the correct
$\log$-corrections for $d=2.$

We define the mass $m_{\varepsilon}\left(  x\right)  ,\;x\in S^{d-1}$, by%
\[
m_{\varepsilon}\left(  x\right)  \overset{\mathrm{def}}{=}-\lim_{k\rightarrow
\infty}\frac{1}{k}\log\lim_{V\uparrow\mathbb{Z}^{d}}\hat{E}_{V,\varepsilon
}\left(  \phi_{0}\phi_{\left[  kx\right]  }\right)  .
\]
The most precise results we have are for the harmonic case:

\begin{theorem}
\label{Th_pinning}

\begin{enumerate}
\item[a)] If $d=2,$ then for small enough $\varepsilon:$%
\[
\left|  \lim_{V\uparrow\mathbb{Z}^{d}}\hat{E}_{V,\varepsilon}^{\mathrm{harm}%
}\left(  \phi_{0}^{2}\right)  -\frac{\left|  \log\varepsilon\right|  }%
{2\pi\sqrt{\det\mathcal{Q}}}\right|  \leq\operatorname*{const}\times
\log\left|  \log\varepsilon\right|
\]

\item[b)] If $d=2,$ then for all $x\in S^{d-1}$ and small enough
$\varepsilon:$%
\[
\operatorname*{const}\times\frac{\sqrt{\varepsilon}}{\left|  \log
\varepsilon\right|  ^{3/4}}\leq m_{\varepsilon}^{\mathrm{harm}}\left(
x\right)  \leq\operatorname*{const}\times\frac{\sqrt{\varepsilon}}{\left|
\log\varepsilon\right|  ^{3/4}}.
\]

\item[c)] If $d\geq3$, then for all $x\in S^{d-1}$ and small enough
$\varepsilon:$%
\[
\operatorname*{const}\times\sqrt{\varepsilon}\leq m_{\varepsilon
}^{\mathrm{harm}}\left(  x\right)  \leq\operatorname*{const}\times
\sqrt{\varepsilon}.
\]
\end{enumerate}

The constants depend on the dimension $d$ and $p$ only.
\end{theorem}

The proof of the results depends on a comparison of the laws of the trapping
configurations with Bernoulli measures. This is particularly delicate in $d=2
$. The following result is the key comparison of the distribution of traps
with Bernoulli measures. We formulate it only in the harmonic case. Somewhat
weaker results are proved in \cite{BoVe} also for the anharmonic situation.

\begin{theorem}
\label{Th_Comparison}Let $\mathcal{A}_{\varepsilon,V}$ be a random subset of
$V$ with $P\left(  \mathcal{A}_{\varepsilon,V}=V\backslash A\right)
=\nu_{V,\varepsilon}\left(  A\right)  .$ Assume $d=2,$ and $U\left(  x\right)
=x^{2}/2.$

\begin{enumerate}
\item[a)] Let $\alpha>0.$ There exists $\varepsilon_{0}>0$ and $C\left(
\alpha\right)  >0$ such that for $\varepsilon\leq\varepsilon_{0}$, any finite
set $V\subset\mathbb{Z}^{d},$ and any $B\subset V$ with $\operatorname*{dist}%
\left(  B,V^{c}\right)  >\varepsilon^{-\alpha},$ one has the estimate%
\[
P\left(  \mathcal{A}_{\varepsilon,V}\cap B=\emptyset\right)  \geq\left(
1-C\left(  \alpha\right)  \frac{\varepsilon}{\sqrt{\left|  \log\varepsilon
\right|  }}\right)  ^{\left|  B\right|  }.
\]

\item[b)] There exist $C>0$ and $\varepsilon_{0}>0$ such that for
$\varepsilon\leq\varepsilon_{0}$, any finite set $V\subset\mathbb{Z}^{d},$ and
all $B\subset V,$ one has%
\[
P\left(  \mathcal{A}_{\varepsilon,V}\cap B=\emptyset\right)  \leq\left(
1-C\frac{\varepsilon}{\sqrt{\left|  \log\varepsilon\right|  }}\right)
^{\left|  B\right|  }%
\]
\end{enumerate}
\end{theorem}

The case of dimension $d\geq3$ is simpler and somewhat better estimates can be
obtained. With the help of the above theorem and the random walk
representation (\ref{RWRepresentation}), a comparison can be made, relating
the quantities in Theorem \ref{Th_pinning} to random trapping problems for
Bernoulli traps. For instance, when investigating the variance, we get%
\[
\hat{E}_{V,\varepsilon}^{\mathrm{harm}}\left(  \phi_{0}^{2}\right)
=E_{\mathrm{traps}}\mathbb{E}_{0}\left(  \sum_{t=0}^{\tau-1}1_{0}\left(
\eta_{t}\right)  \right)  =E_{\mathrm{traps}}\sum_{t=0}^{\infty}p_{t}\left(
0\right)  \mathbb{P}_{0,0}^{\left(  t\right)  }\left(  \mathcal{A\cap}%
\eta_{\left[  0,t\right]  }=\emptyset\right)  ,
\]
where $\mathcal{A}$ is the random set of points with traps, as introduced
above, $\tau$ is the first entrance time into this trapping set and
$\eta_{\left[  0,t\right]  }$ is the set of points visited by the random walk
between time $0$ and $t$. $\mathbb{P}_{0,0}^{\left(  t\right)  }$ refers to a
random walk bridge from $0$ to $0$ in time $t,$ and $p_{t}\left(  x\right)
,\;x\in\mathbb{Z}^{d}$ are the transition probabilities of the random walk.
With the help of Theorem \ref{Th_Comparison}, the right hand side can be
estimated in terms of a Bernoulli trapping problem. If the traps are Bernoulli
on $\mathbb{Z}^{d},$ with probability $\rho$ that a trap is present at a given
site, then (in the $V\uparrow\mathbb{Z}^{d}$ limit)%
\[
E_{\mathrm{traps}}\sum_{t=0}^{\infty}p_{t}\left(  0\right)  \mathbb{P}%
_{0,0}^{\left(  t\right)  }\left(  \mathcal{A\cap}\eta_{\left[  0,t\right]
}=\emptyset\right)  =\sum_{t=0}^{\infty}p_{t}\left(  0\right)  \mathbb{E}%
_{0,0}^{\left(  t\right)  }\exp\left[  \left|  \eta_{\left[  0,t\right]
}\right|  \log\left(  1-\rho\right)  \right]  .
\]
There are classical results about the right hand side, due to Donsker and
Varadhan \cite{DoVa}, Sznitman \cite{Sz}, and most recently in
\cite{vdBeBodHo} investigating such questions. The classical Donsker-Varadhan
result is not sharp enough to prove the results of Theorem \ref{Th_pinning},
but a modification of the arguments in \cite{vdBeBodHo} is exactly what is
needed. The following result is a discrete but somewhat weaker version of one
of the main results in \cite{vdBeBodHo}.

\begin{proposition}
Assume $d=2.$ There exists a function $\mathbb{R}^{+}\ni a\rightarrow r\left(
a\right)  \in\mathbb{R}^{+}$, satisfying $\lim_{a\rightarrow0}r\left(
a\right)  =\infty,$ such that%
\[
\mathbb{P}_{0,0}^{\left(  t\right)  }\left(  \left|  \eta_{\left[  0,t\right]
}\right|  \leq a\frac{t}{\log t}\right)  \leq t^{-r\left(  a\right)  }%
\]
for large enough $t.$
\end{proposition}

(In \cite{vdBeBodHo}, a variational formula for $r\left(  a\right)  $ is
given, in the continuous Wiener sausage case.) This proposition and Theorem
\ref{Th_Comparison} lead to the appropriate variance estimates in Theorem
\ref{Th_pinning} a).

For the anharmonic case, the results are less precise, but we still get the
correct leading order dependence of the variance on the pinning parameter
$\varepsilon.$ Assume that there is a $C>0$ such that%
\[
1/C\leq U^{\prime\prime}\left(  x\right)  \leq C,\;\forall x.
\]
Under this condition we have the following result:

\begin{theorem}
Assume $d=2.$ There exists a constant $D$, depending on $p,$ such that%
\[
\frac{1}{D}\left|  \log\varepsilon\right|  \leq\sup_{V}\hat{E}_{V,\varepsilon
}\left(  \phi_{0}^{2}\right)  \leq D\left|  \log\varepsilon\right|  .
\]
\end{theorem}

The upper bound is in \cite{DeVe} and \cite{IoVe}, and the lower bound is in
\cite{BoVe}.

\section{Entropic repulsion and the wetting transition\label{Chap_EntropicRep}}

\vskip-5mm \hspace{5mm}

In view of the example of Fisher discussed shortly in Section \ref{Chap_Intro}
it is natural to ask similar question for higher-dimensional interfaces. The
first task is to investigate the effect of a wall on the random surface
without the presence of a pinning effect. There are different ways to take the
presence of a wall into account. We have mainly worked with a ``hard wall'',
i.e. where the measure is simply conditioned on the event $\Omega_{V}%
^{+}\overset{\mathrm{def}}{=}\left\{  \phi:\phi_{x}\geq0\;\forall x\in
V\right\}  .$ There are other possibilities, for instance by introducing a
``soft wall''. This means that the Hamiltonian (\ref{Hamiltonian}) is changed
by adding $\sum_{x\in V}f\left(  \phi_{x}\right)  ,$ where
$f:\mathbb{R\rightarrow R}$ satisfies $\lim_{x\rightarrow\infty}f\left(
x\right)  =0,$ $\lim_{x\rightarrow-\infty}f\left(  x\right)  =\infty.$ We will
only work with a hard wall here, and consider the conditional law for the
random field $P\left(  \cdot|\Omega_{V}^{+}\right)  $.

What is the effect of the presence of the wall on the surface? The crucial
point is that the surface has local fluctuations, which push the interface
away from the wall. On the other hand, the long-range correlations give the
surface a certain global stiffness. In order to understand what is going on,
consider first the case where there are no such long-range correlations, in
the extreme case, where the $\phi_{x}$ are just i.i.d. random variables. In
that case, evidently nothing interesting happens: The variables are
individually conditioned to stay positive. In particular, $E\left(  \phi
_{x}|\Omega_{V}^{+}\right)  $ stays bounded for $V\uparrow\mathbb{Z}^{d}.$
This picture remains the same for fields with rapidly decaying correlations.
However, gradient fields behave differently, and so do interfaces in more
realistic statistical physics models. As the surface has some global
stiffness, the energetically best way for the surface to leave some room for
the local fluctuations is to move away from the wall in some global sense.
This effect is called ``entropic repulsion'' and is well known in the physics literature.

The first mathematically rigorous treatment of entropic repulsion appeared in
the paper by Bricmont, Fr\"{o}hlich and El Mellouki \cite{BrFrMe}. In a series
of papers \cite{BoDeZe1}, \cite{De}, \cite{DeGi1}, and \cite{BoDeGi}, sharp
quantitative results have been derived, the most accurate ones for the
harmonic case.

In most of these and related questions, the two-dimensional case is the most
difficult but also the most interesting one. In fact, interfaces in the ``real
world'' are mostly two-dimensional.

We first present the results for $d\geq3.$ For gradient non-Gaussian models,
some results in the same spirit have been obtained in \cite{DeGi1}, but they
are not as precise as the ones obtained in the Gaussian model. The case where
one starts with the field $P_{\infty}$ (which exists for $d\geq3)$ is somewhat
easier than the field on the finite box $V_{n}$ with zero boundary condition.
In the latter case, there are some boundary effects complicating the situation
without changing it substantially. This is investigated in \cite{De}. Despite
the fact that we consider $P_{\infty}$, we consider the wall only on a finite
box, i.e., we consider $P_{\infty}\left(  \cdot|\Omega_{V_{n}}^{+}\right)  ,$
and we are interested in what happens as $n\rightarrow\infty.$ We usually
write $\Omega_{n}^{+}$ for $\Omega_{V_{n}}^{+}.$ Our first task is to get
information about $P_{\infty}\left(  \Omega_{n}^{+}\right)  .$ The following
results are proved only for the case of nearest neighbor interactions, i.e.
when $p\left(  x\right)  =1/2d$ for $\left|  x\right|  =1.$

\begin{theorem}
\label{ThEntrRep}Let $d\geq3$. Then

\begin{enumerate}
\item[a)]
\[
P_{\infty}^{\mathrm{harm}}\left(  \Omega_{n}^{+}\right)  =\exp\left[
-2\Gamma(0)\mathrm{cap}\left(  V\right)  n^{d-2}\log n\left(  1+o(1)\right)
\right]  ,
\]
where $V=\left[  -1,1\right]  ^{d},$ $\mathrm{cap}(A)$\textrm{\ }denotes the
Newtonian capacity of $A$%
\[
\mathrm{cap}(A)\overset{\mathrm{def}}{=}\inf\left\{  \left\|  \nabla
f\right\|  ^{2}:f\geq1_{A}\right\}  ,
\]
and $\Gamma(0)=\gamma_{\infty}(0,0)$ is the variance of $\phi_{0}$ under
$P_{\infty}^{\mathrm{harm}}.$

\item[b)]
\[
E_{\infty}^{\mathrm{harm}}\left(  \phi_{0}|\Omega_{n}^{+}\right)
=2\sqrt{\Gamma(0)\log n}(1+o(1)).
\]

\item[c)]
\[
\mathcal{L}_{P_{\infty}^{\mathrm{harm}}(\cdot|\Omega_{n}^{+})}\left(  \left(
\phi_{x}-E_{\infty}(\phi_{x}|\Omega_{n}^{+})\right)  _{x\in\mathbb{Z}^{d}%
}\right)  \rightarrow P_{\infty}^{\mathrm{harm}}\text{\textrm{\ weakly,}}%
\]
as $n\rightarrow\infty,$ where $\mathcal{L}_{P_{\infty}^{\mathrm{harm}}%
(\cdot|\Omega_{n}^{+})}$ denotes the law of the field under the conditioned measure.
\end{enumerate}
\end{theorem}

Part b) gives the exact rate at which the random surface escapes to infinity,
while part c) states that the effect of the entropic repulsion essentially is
only this shifting: after subtraction of the shift by the expectation, the
surface looks as it does without the wall. However, there is some subtlety in
this picture. From the theorem in particular part c), one might conclude that
$\lim_{n\rightarrow\infty}P_{\infty}^{\mathrm{harm}}\theta_{2\sqrt
{\Gamma(0)\log n}}^{-1}\left(  \Omega_{n}^{+}\right)  =1,$ where $\theta
_{a}:\mathbb{R}^{\mathbb{Z}^{d}}\rightarrow\mathbb{R}^{\mathbb{Z}^{d}}$ is the
shift mapping $\theta_{a}\left(  \left(  \phi_{x}\right)  _{x\in\mathbb{Z}%
^{d}}\right)  =\left(  \phi_{x}+a\right)  _{x\in\mathbb{Z}^{d}}.$ But this is
not the case. In fact $P_{\infty}^{\mathrm{harm}}\theta_{2\sqrt{\Gamma(0)\log
n}}^{-1}\left(  \Omega_{n}^{+}\right)  $ converges rapidly to $0.$ As part c)
states only the weak convergence, this is no contradiction. Parts a) and b) of
Theorem \ref{ThEntrRep} had been proved in \cite{BoDeZe1}, part c) in
\cite{DeGi1}.

We come now to the two-dimensional case which is considerably more delicate
than the higher dimensional one. We again consider only the harmonic case. We
write $P_{n}$ for $P_{V_{n}}.$ If the lattice is two-dimensional, a
thermodynamic limit of the measures $P_{n}$ does not exist as the variance
blows up. $P_{n}^{\mathrm{harm}}(\Omega_{n}^{+})$ is of order $\exp\left[
-cn\right]  ,$ as has been shown in \cite{De}. As remarked above, this is
mainly a boundary effect and is not really relevant for the phenomenon of the
entropic repulsion. To copy somehow the procedure of the case $d\geq3$, we
consider a subset $D\subset V=[-1,1]^{2}$ which has a nice boundary and a
positive distance from the boundary of $V.$ To be specific, just think of
taking $D\overset{\mathrm{def}}{=}\lambda V$ for some $\lambda<1.$ Then let
$D_{n}\overset{\mathrm{def}}{=}nD\cap\mathbb{Z}^{2}$ and $\Omega_{D_{n}}%
^{+}\overset{\mathrm{def}}{=}\left\{  \phi_{x}\geq0,\,x\in D_{n}\right\}  .$
In contrast to $P_{n}(\Omega_{n}^{+}),$ $P_{n}\left(  \Omega_{D_{n}}%
^{+}\right)  $ decays much slower, but still faster than any polynomial rate.
In \cite{BoDeGi} we proved the following result:

\begin{theorem}
\label{ThEntrRep2d}Assume $d=2$ and let $g\overset{\mathrm{def}}{=}1/2\pi.$

\begin{enumerate}
\item[a)]
\[
\lim_{n\rightarrow\infty}\frac{1}{(\log n)^{2}}\log P_{n}^{\mathrm{harm}%
}(\Omega_{D_{n}}^{+})=-2g\mathrm{cap}_{V}(D),
\]
where $\mathrm{cap}_{V}(D)$ is the relative capacity of $D$ with respect to
$V$:%
\[
\mathrm{cap}_{V}(D)\overset{\mathrm{def}}{=}\inf\left\{  \left\|  \nabla
f\right\|  _{2}^{2}:f\in H_{0}^{1}(V),\,f\geq1\,\mathrm{on\,}D\right\}  .
\]
Here, $H_{0}^{1}(V)$ is the Sobolev space of (weakly) differentiable functions
$f$ with square integrable gradient and $f|_{\partial V}=0.$

\item[b)] For any $\varepsilon>0$%
\[
\lim_{n\rightarrow\infty}\sup_{x\in D_{n}}P_{n}^{\mathrm{harm}}\left(  \left.
\left|  \phi_{x}-2\sqrt{g}\log n\right|  \geq\varepsilon\log n\right|
\Omega_{D_{n}}^{+}\right)  =0.
\]
\end{enumerate}
\end{theorem}

This corresponds to parts a) and b) of Theorem \ref{ThEntrRep}. Part c) does
not make sense here as $P_{\infty}^{\mathrm{harm}}$ does not exist. Remark
that under the unconditional law $P_{n}^{\mathrm{harm}},$ $\left|  \phi
_{x}\right|  $ is typically of order $\sqrt{\log n}$ in the bulk.

Roughly speaking, the delicacy in the two-dimensional case is coming from the
fact that the relevant ``spikes'' responsible for the repulsion are thicker
than in the higher dimensional case, where essentially just very local spikes
are responsible for the effect. This makes necessary to apply a multiscale
analysis separating the scales of the spikes.

It is well-known that the two-dimensional harmonic field has much similarity
with a hierarchical field defined in the following way: We call a sequence
$\alpha=\alpha_{1}\alpha_{2}\ldots\alpha_{m}$, $\alpha_{i}\in\{0,1\}$ a
\textit{binary string.} $\ell(\alpha)=m$ is the length. $\emptyset$ is the
empty string of length $0.$ We write $T$ for the set of all such strings of
finite length, and $T_{m}\subset T$ for the set of strings of length $m.$ If
$\alpha\in T_{m},\,0\leq k\leq m,$ we write $[\alpha]_{k}$ for the truncation
at level $k:$%
\[
\left[  \alpha_{1}\alpha_{2}\ldots\alpha_{m}\right]  _{k}\overset
{\mathrm{def}}{=}\alpha_{1}\alpha_{2}\ldots\alpha_{k}.
\]
If $\alpha,$ $\beta\in T_{m}$ we define the hierarchical distance
\[
d_{H}(\alpha,\beta)\overset{\mathrm{def}}{=}m-\max\left\{  k\leq
m:[\alpha]_{k}=[\beta]_{k}\right\}.
\]
We consider the following family $(X_{\alpha})_{\alpha\in T_{m}}$ of centered
Gaussian random variables by%
\begin{equation}
\mathrm{cov}\left(  X_{\alpha},X_{\beta}\right)  =\gamma\left(  m-d_{H}%
(\alpha,\beta)\right)  ,\label{covhier}%
\end{equation}
with a parameter $\gamma>0.$ We argue now that there is much similarity
between the two dimension harmonic field $(\phi_{x})_{x\in D_{n}}$ and the
field $\left(  X_{\alpha}\right)  _{\alpha\in T_{m}}.$ To see this, we first
match the number of variables, i.e. put $2^{m}=\left|  D_{n}\right|  .$ As
$\left|  D_{n}\right|  $ is of \ order $n^{2},$ this just means that
$m\sim2\log n/\log2.$ Then we should also match the variances, i.e. take
$\gamma=g/2\log2.$ For the free field $(\phi_{x})$, it is known that
$\mathrm{cov}(\phi_{x},\phi_{y})$ behaves like $g\left(  \log n\right)
/\log\left|  x-y\right|  $, if $x,y$ are not too close to the boundary. This
follows from the random walk representation. Comparing this with
(\ref{covhier}), we see that for any number $s\in(0,g)$%
\begin{equation}
\#\left\{  y\in D_{n}:\mathrm{cov}\left(  \phi_{x},\phi_{y}\right)  \leq s\log
n\right\}  \sim\#\left\{  \beta\in T_{m}:\mathrm{cov}\left(  X_{\alpha
},X_{\beta}\right)  \leq s\log n\right\}  \label{free=hier}%
\end{equation}
in first order, for any $x\in D_{n},$ $\alpha\in T_{m}.$ Therefore, the two
fields have roughly the same covariance structure. The hierarchical field is
much simpler and is very well investigated (see e.g. \cite{Bi}, \cite{Bramson}%
, \cite{DerridaSpohn}), and the entropic repulsion is much easier to discuss
than for the harmonic field. The approach to prove Theorem \ref{ThEntrRep2d}
consist in introducing a hierarchical structure in the $\left(  \phi
_{x}\right)  $-field with the help of successive conditionings on a hierarchy
of scales, and then adapt the methods from the purely hierarchical case.

We come now back to the question of a wetting transition, as discussed in the
one-dimensional case by Michel Fisher \cite{Fisher}. One is interested in the
behavior of $\hat{P}_{V,\varepsilon}\left(  \cdot\mid\Omega_{V}^{+}\right)  $
for large $V,$ where $\hat{P}_{V,\varepsilon}$ is the pinned measure
introduced in (\ref{Def_PV-pinned}). Unfortunately, we are not able to
describe this path measure. The simplest way to discuss the wetting transition
is in terms of free energy considerations. For this we expand $\hat
{P}_{V,\varepsilon}\left(  \Omega_{V}^{+}\right)  $ (see (\ref{Expand_Pinning}%
)):%
\[
\hat{P}_{V,\varepsilon}\left(  \Omega_{V}^{+}\right)  =\sum_{A\subset
V}\varepsilon^{\left|  V\backslash A\right|  }\frac{Z_{A}}{\hat{Z}%
_{V,\varepsilon}}P_{A}\left(  \Omega_{V}^{+}\right)  .
\]
It is plausible, that pinning ``wins'' over entropic repulsion, if this sum is
much larger than the contribution to the sum coming from subsets $A$ having
essentially no pinning sites, i.e. $A\approx V.$ It is therefore natural to
consider the quantity%
\[
p_{+}\left(  \varepsilon\right)  \overset{\mathrm{def}}{=}\lim_{V\uparrow
\mathbb{Z}^{d}}\frac{1}{\left|  V\right|  }\log\frac{\hat{Z}_{V,\varepsilon
}\hat{P}_{V,\varepsilon}\left(  \Omega_{V}^{+}\right)  }{Z_{V}P_{V}\left(
\Omega_{V}^{+}\right)  }=\lim_{V\uparrow\mathbb{Z}^{d}}\frac{1}{\left|
V\right|  }\log\frac{\hat{Z}_{V,\varepsilon}\hat{P}_{V,\varepsilon}\left(
\Omega_{V}^{+}\right)  }{Z_{V}}.
\]
The limit is easily seen to exist. It is also not difficult to see that
$p_{+}\left(  \varepsilon\right)  >0$ for large enough $\varepsilon>0$, and in
any dimension (see \cite{BoIo}). Similar to the discrete random walk case in
\cite{Fisher}, the Gaussian model has a wetting transition, too, for $d=1$:
There exists an $\varepsilon_{\mathrm{crit}}>0,$ such that $p_{+}\left(
\varepsilon\right)  =0$ for $\varepsilon<\varepsilon_{\mathrm{crit}}.$ This is
easy to see for $d=1.$ For the harmonic model, there is remarkably no such
transition for $d\geq3,$ but for $d=2$ there is a wetting transition.

\begin{theorem}
{\rm \cite{BoDeZe2}} For $d\geq3,$ $p_{+}^{\mathrm{harm}}\left(
\varepsilon\right)
>0$ for all $\varepsilon>0.$
\end{theorem}

\begin{theorem}
{\rm \cite{CaVe}} For $d=2,$ there exists $\varepsilon_{\mathrm{crit}}%
^{\mathrm{harm}}>0,$ such that $p_{+}^{\mathrm{harm}}\left(  \varepsilon
\right)  =0$ for $\varepsilon<\varepsilon_{\mathrm{crit}}^{\mathrm{harm}}.$
\end{theorem}

Remarkably, too, Caputo and Velenik have proved that such a wetting transition
exists for $d\geq3$ for some non-harmonic models, e.g. for $U\left(  x\right)
=\left|  x\right|  .$

There are many open questions concerning this wetting transition, which is
very poorly understood (mathematically). For instance, the methods discussed
in Section \ref{Chap_Pinning} do not apply, and we are not able to prove that
in the pinning dominated region $p_{+}\left(  \varepsilon\right)  >0,$ the
measure is pathwise localized, i.e. that%
\[
\sup_{V}\sup_{x\in V}\hat{P}_{V,\varepsilon}\left(  \phi_{x}^{2}\mid\Omega
_{V}^{+}\right)  <\infty,
\]
which certainly should be expected. To discuss the nature of the transition
(first order or second order?) is probably even much more delicate.

\section{Localization-delocalization transitions for one-dimensional
copolymers\label{Chap_Copolymer}}

\vskip-5mm \hspace{5mm}

We stick here to the standard simple random walk case where $P_{n}$ simply is
the uniform distribution on the set of paths $\phi_{0}=0,\phi_{1},\ldots
,\phi_{n}\in\mathbb{Z},$ satisfying $\left|  \phi_{i}-\phi_{i-1}\right|
=1,\;1\leq i\leq n.$ There is not much difference when considering more
general random walks, or the tied-down situation, but most of the published
results are for the simple random walk. An interesting case of a mixed
attractive-repulsive interaction is given in the following way. Regard the
above random walk as a (very simplified) model of a polymer chain imbedded in
two liquids, say water and oil. The water is at the bottom, say at points
$\left(  i,j\right)  \in\mathbb{N\times Z},$ $j\leq0,$ and the oil above at
$j>0.$ The polymer chain is attached with one end at the interface between the
two liquids, and interacts with them in the following way: To each ``node''
$\left(  i,\phi_{i}\right)  $ of the polymer chain, we attach a value
$\sigma_{i}\in\mathbb{R}$ which is $<0$ if the node is water-repellent, and
$>0$ if it is oil-repellent. The overall effect is described by the
Hamiltonian%
\[
H_{n,\sigma}\left(  \phi\right)  \overset{\mathrm{def}}{=}\sum_{i=1}^{n}%
\sigma_{i}\operatorname*{sign}\left(  \phi_{i}\right)  ,
\]
where we put $\operatorname*{sign}\left(  0\right)  \overset{\mathrm{def}}%
{=}0.$ With this Hamiltonian, we define the $\sigma$-dependant path measure%
\[
P_{n,\beta,\sigma}\left(  \phi\right)  \overset{\mathrm{def}}{=}\frac
{1}{Z_{n,\beta,\sigma}}\exp\left[  -\beta H_{n,\sigma}\left(  \phi\right)
\right]  ,
\]
where $\beta>0$ is a parameter governing the strength of the interaction. We
assume that the $\sigma_{i}$ change sign either in a periodic way or randomly.
There may be two competing effects. The polymer chain may try to follow the
preferences described by the $\sigma$'s as closely as possible in which case
the path evidently would have to stay close to the oil-water-interface and
gets localized. On the other hand, this strategy may be entropically too
costly, in particular if there is no balance between oil-repellence and
water-repellence. We will always assume that%
\[
h\overset{\mathrm{def}}{=}\lim_{n\rightarrow\infty}\frac{1}{n}\sum_{j=1}%
^{n}\sigma_{j},
\]
exists, and we assume it to be $\geq0.$ (The case $h\leq0$ can be treated
symmetrically). It turns out that typically, there is a non-trivial curve in
the $\left(  \beta,h\right)  $-plane which separates the localized from the
delocalized region. This phase separation line is quite model dependent, but
the behavior near $\left(  0,0\right)  $ appears to be much more universal but
it is completely different depending whether the $\sigma_{i}$ are random or periodic.

The first rigorous results in this direction had been obtained by Sinai
\cite{Sinai} who proved the following result in the balanced case (i.e.
$h=0).$ Let $\mathbb{P}$ be the symmetric Bernoulli-measure on $\left\{
-1,1\right\}  ^{\mathbb{N}}.$

\begin{theorem}
Let $\beta>0.$ There exist constants $C$ and $\rho\left(  \beta\right)  >0,$
and for $\mathbb{P}$-almost all $\sigma=\left(  \sigma_{i}\right)  _{i\geq1},$
there exists a sequence $\left(  R_{n}\left(  \sigma\right)  \right)
_{n\in\mathbb{N}}$ of natural numbers such that%
\[
P_{n,\beta,\sigma}\left(  \left|  \phi_{n}\right|  \geq r\right)  \leq
C\exp\left[  -\rho\left(  \beta\right)  r\right]
\]
for $r\geq R_{n}\left(  \sigma\right)  .$ The sequence $\left(  R_{n}\right)
$ is stochastically bounded, i.e.%
\[
\lim_{m\rightarrow\infty}\sup_{n}\mathbb{P}\left(  R_{n}\geq m\right)  =0.
\]
\end{theorem}

In a paper with Frank den Hollander \cite{BodHo} we proved that there is a
localization-delocalization transition in the random non-balanced case. This
transition is discussed in this paper in terms of the free energy. To describe
the results, let $\sigma_{i}=\pm1+h$ with probabilities $1/2,$ and
independently, $h\geq0.$ One strategy of the path could be just to stay on the
negative side all the time, i.e. $\phi_{i}<0$ for all $i\leq n.$ This leads to
a trivial lower bound of the free energy%
\[
f\left(  \beta,h\right)  \overset{\mathrm{def}}{=}\lim_{n\rightarrow\infty
}\frac{1}{n}\log Z_{n,\beta,\sigma}%
\]
which is easily seen to exist, and is non-random:%
\begin{align*}
Z_{n,\beta,\sigma}  & \geq E_{n}\left(  \exp\left[  -\beta H_{n,\sigma}\left(
\phi\right)  \right]  1_{\left\{  \phi_{i}<0,\;\forall i\leq n\right\}
}\right)  \\
& =\exp\left[  \beta\sum\nolimits_{i=1}^{n}\sigma_{i}\right]  P\left(
\phi_{i}<0,\;\forall i\leq n\right)  .
\end{align*}
From this we get%
\[
f\left(  \beta,h\right)  \geq\beta h.
\]
It is quite plausible that localization dominates in the case where there is a
strict inequality, and that delocalization holds if $f\left(  \beta,h\right)
=\beta h.$

\begin{theorem}
There exists a positive, continuous, and increasing function $\beta\rightarrow
h^{\ast}\left(  \beta\right)  $ such that%
\begin{align}
f\left(  \beta,h\right)   & >\beta h\;\mathrm{for\;}0\leq h<h^{\ast}\left(
\beta\right)  ,\label{Copolymer_loc}\\
f\left(  \beta,h\right)   & =\beta h\;\mathrm{for\;}h>h^{\ast}\left(
\beta\right)  .\label{Copolymer_deloc}%
\end{align}
The function $\beta\rightarrow h^{\ast}\left(  \beta\right)  $ has a positive
tangent at $\beta=0.$
\end{theorem}

The phase separating function $h^{\ast}$ is certainly very much model
dependent, but we expect that the tangent at $0$ is model independent, and
would be the same for any random law of the $\sigma$-sequence which has
variance $1$ and a expectation $h$, and has exponentially decaying tails, but
this is not proved in \cite{BodHo}. In physics literature, there are
non-rigorous arguments claiming that the tangent is $1,$ but we neither have
been able to prove or disprove it, yet. We prove that the tangent at $0$ can
be described in terms of a phase separation line for a continuous model, where
the random walk is replaced by a Brownian motion, and the random environment
$\sigma$ is replaced by (biased) white noise. In this case, the phase
separation line is a straight line, and we prove that this line is the tangent
at $0$ of our model. It should be remarked that the $\left(  \beta,h\right)
\approx\left(  0,0\right)  $ situation, cannot be handled by simple
perturbation techniques.

A natural question is if in the localized region $f\left(  \beta,h\right)
>\beta h$ the path measure is really localized in the sense described in the
paper of Sinai. This is indeed the case and has been proved by Biskup and den
Hollander \cite{BidHo}. One might also wonder if in the localized region
$f\left(  \beta,h\right)  =\beta h$ or at least in the interior of it, the
path measure is really delocalized, which should mean, that it converges,
after Brownian rescaling, to the limit of a random walk conditioned to stay
negative, which is the negative of the Brownian meander. This seems to be a
rather difficult question and has not been answered, yet.

The positive tangent is essentially tied to the randomness of the sequence.
For the periodic case, the situation is different, as has recently been proved
in \cite{BoGi}:

\begin{theorem}
Let $\sigma_{i}=\omega_{i}+h,$ where $\omega_{i}\in\left\{  -1,1\right\}  $ is
periodic, i.e. such that there exists $T$ with $\omega_{i+2T}=\omega_{i}$ for
all $i,$ and $\sum_{i=1}^{2T}\omega_{i}=0.$ Then there is a function $h^{\ast
}$ such that (\ref{Copolymer_loc}) and (\ref{Copolymer_deloc}) hold. In this
case%
\[
C=\lim_{\beta\rightarrow0}\frac{h^{\ast}\left(  \beta\right)  }{\beta^{3}}%
\]
exists and is positive.
\end{theorem}

In this paper an expression for $C$ in terms of a variational problem is
derived, where the exact nature of the periodic sequence enters.

\label{lastpage}


\begin{thebibliography}{99}
                                                                                              %
\bibitem {vdBeBodHo}van den Berg, M., Bolthausen, E. and den Hollander, F.:
\textit{Moderate deviations for the Wiener sausage}. Annals of Mathematics \textbf{153}(2001), 355--406.

\bibitem {Bi}Biggins, J.D.: \textit{Chernoff's theorem in the branching random
walk}. J. Appl. Prob. \textbf{14}(1977), 630--636.

\bibitem {BidHo}Biskup, M. and den Hollander, F.: \textit{A heteropolymer near
a linear interface}. Ann. Appl. Probab. \textbf{9}(1999),
668--687.

\bibitem {BoDeZe1}Bolthausen, E., Deuschel, J. D., and Zeitouni, O.:
\textit{Entropic repulsion for the lattice free field}, Comm.
Math. Phys. \textbf{170}(1995), 417--443.

\bibitem {BoDeZe2}Bolthausen, E., Deuschel, J.D., and Zeitouni, O.
\textit{Absence of a wetting transition for lattice free fields in
dimensions three and larger}. J. Math. Phys. \textbf{41}(2000),
1211--1223.

\bibitem {BoDeGi}Bolthausen, E., Deuschel, J.D., and Giacomin, G.:
\textit{Entropic repulsion for the two-dimensional lattice free
field}. Annals of Probability \textbf{29}(2001), 1670--1692.

\bibitem {BoBr}Bolthausen, E., and Brydges, D.: \textit{Gaussian surface
pinned by a weak potential}. IMS Lecture Notes Vol.
\textbf{36}(2001), 134--149.

\bibitem {BodHo}Bolthausen, E. and den Hollander, F.: \textit{Localization
transition for a polymer near an interface}. Ann. Probability
\textbf{25}(1997), 1334--1365.

\bibitem {BoIo}Bolthausen, E. and Ioffe, D.: \textit{Harmonic crystal on the
wall: a microscopic approach}, Comm. Math. Phys.
\textbf{187}(1997), 523--566.

\bibitem {BoVe}Bolthausen, E. and Velenik, Y.: \textit{Critical behavior of
the massless free field at the depinning transition. }Comm. Math.
Phys. \textbf{233}(2001), 161--203.

\bibitem {BoGi}Bolthausen, E. and Giacomin, G.: \textit{On the critical
delocalization-localization line for periodic copolymers at
interfaces. }Preprint.

\bibitem {Bramson}Bramson, M.: \textit{Maximal displacement of branching
Brownian motion.} Comm. Pure Appl. Math. \textbf{31}(1978),
531--581.

\bibitem {Brezin}Brezin, E., Halperin, and Leibler, S.: \textit{Critical
wetting in three dimensions}. Phys. Rev. Lett. \textbf{50}(1983),
1387.

\bibitem {BrFrMe}Bricmont, J., Fr\"{o}hlich, J., and El Mellouki, A.:
\textit{Random surfaces in statistical mechanics: Roughening,
rounding, wetting}. J. Stat. Phys. \textbf{42}(1986), 743.

\bibitem {BrFrSp}Brydges, D.C., Fr\"{o}hlich, J., and Spencer, T.: \textit{The
random walk representation of classical spin systems and
correlation inequalities. }Commun. Math. Phys., \textbf{83}(1982),
123--150.

\bibitem {CaVe}Caputo, P., and Velenik, I.: \textit{A note on wetting
transition for gradient fields}. Stoch. Proc. Appl.
\textbf{87}(2000), 107--113.

\bibitem {De}Deuschel, J.D.: \textit{Entropic repulsion of the lattice free
field.} II. The $0$-boundary case. Comm. Math. Phys.
\textbf{181}(1996), 647--665.

\bibitem {DeGi1}Deuschel, J. D., and Giacomin G.: \textit{Entropic repulsion
for the free field: pathwise characterization in }$d\geq3$, Comm.
Math. Phys. \textbf{206}(1999), 447--462.

\bibitem {DeGi2}Deuschel, J.D., and Giacomin, G.: \textit{Entropic repulsion
for massless fields}, Stoch. Process. Appl. \textbf{89}(2000),
333--354.

\bibitem {DeGiIo}Deuschel, J.D., Giacomin, G., and Ioffe, D.: \textit{Large
deviations and concentration properties for }$\nabla\phi$\textit{\
interface models}. Prob. Theory Rel. Fields \textbf{117}(2000),
49--111.

\bibitem {DeVe}Deuschel, J.D., and Velenik, Y.: \textit{Non-Gaussian surface
pinned by a weak potential}. Prob. Theory Rel. Fields
\textbf{116}(2000), 359--377.

\bibitem {DerridaSpohn}Derrida, B. and Spohn, H.: \textit{Polymers on
disordered trees, spin glasses, and travelling waves, }J. Stat.
Phys. \textbf{51}(1988), 817--840.

\bibitem {DoVa}Donsker, M. and Varadhan, S.R.S.: \textit{On the number of
distinct sites visited by a random walk.} Comm. Pure Appl. Math.
\textbf{32}(1979), 721--747.

\bibitem {DuMaRiRo}Dunlop, F., Magnen, J., Rivasseau, V., and Roche, Ph.:
\textit{Pinning of an interface by a weak potential.} J. Statist.
Phys. \textbf{66}(1992), 71--98.

\bibitem {Fisher}Fisher, M.: \textit{Walks, walls, wetting and melting.} J.
Stat. Phys. \textbf{34}(1984), 667--729.

\bibitem {HeSj}Helffer, B., and Sj\"{o}strand, J.: \textit{On the correlation
for Kac-like models in the convex case.} J. Statist. Phys.
\textbf{74}(1994), 349--409.

\bibitem {IoVe}Ioffe, D., and Velenik, I.: \textit{A note on the decay of
correlations under }$\delta$-\textit{pinning. }Prob. Theory Rel.
Fields \textbf{116}(2000), 379--389.

\bibitem {Jap}Isozaki, Y. and Yoshida, N.: \textit{Weakly pinned random walk
on the wall: pathwise descriptions of the phase transition.
}Stoch. Proc. Appl. \textbf{96}(2001), 261--284.

\bibitem {Sinai}Sinai, Ya. G.: \textit{A random walk with a random potential.}
Theory Probab. Appl. \textbf{38}(1993), 382--385.

\bibitem {Sz}Sznitman, A.-S.: \textit{Brownian Motion, Obstacles, and Random
Media. }Springer, Heidelberg 1998.
\end{thebibliography}
\end{document}